\documentclass[a4paper]{article}
\usepackage{amsmath, amssymb, eucal, amscd}

\newtheorem{thm}{Theorem}
\newtheorem{lem}[thm]{Lemma}
\newtheorem{eg}[thm]{Example}
\newtheorem{prop}[thm]{Proposition}
\newtheorem{cor}[thm]{Corollary}
\newtheorem{defn}[thm]{Definition}
\newtheorem{rem}[thm]{Remark}

\newcommand{\smnoind}{{\smallskip\noindent}}
\newcommand{\tar}{{\textbf t}}
\newcommand{\s}{{\textbf s}}
\newcommand{\id}{{\rm id}}

\newcommand{\ti}{\tilde}

\newenvironment{prf}{{\noindent \textbf{Proof:} }}{\hfill $\Box$\medskip}

\begin{document}

\title{Some remarks on groupoids and small categories}
\author{Chi-Keung Ng}
\date{}

\maketitle

\begin{abstract}
This unpublished note contains some materials taken from my old study note on groupoids and small categories (\cite{Ng-note-gpd}). 
It contains a proof for the fact that a groupoid is any group bundle over an equivalence relation. 
Moreover, the action of a category $G$ on a category $H$ as well as the resulting semi-direct product category $H\times_\alpha G$ will be defined (when either $G$ is a groupoid or $H^{(0)} = G^{(0)}$). 
If both $G$ and $H$ are groupoids, then $H\times_\alpha G$ is also a groupoid. 
The reason of producing this note is for people who want to check some details in a recent work of Li (\cite{li}).
\end{abstract}

\bigskip

We consider a small category as a set $H$ (of morphisms) together with the source and the target maps $\s,\tar: H \rightarrow H^{(0)}\subseteq H$ (i.e., we identify elements in $H^{(0)}$ with their identity morhpisms) as well as the composition $(h,h') \mapsto hh'$ as a map from $H^{(2)} := \{(h,h')\in H\times H: \s(h) = \tar(h') \}$ to $H$. 
We recall that a groupoid is a small category such that every morphism has an inverse. 

\medskip

\begin{defn}
Let $X$ be a set and $R$ be an equivalence relation on $X$.
Suppose that $G_\xi$ is a group for any $\xi\in X/R$.
Then $(X, R, \{G_\xi\}_{\xi\in X/R})$ (or simply $\{G_\xi\}_{\xi\in X/R}$)
is called a \emph{group bundle over the equivalence classes of $R$}.
\end{defn}

\medskip

Let $\{G_\xi\}_{\xi\in X/R}$ be a group bundle over the equivalence classes of an equivalence relation $R$ on $X$ and let
$$G\ :=\ \{(x,g,y): \xi\in X/R;\ x,y\in \xi;\ g\in G_\xi\}$$
and $\tar,\s: G\rightarrow X$ are defined by $\tar(x,g,y) = x$ and $\s(x,g,y) = y$.
Moreover, we have $i: X \rightarrow G$ given by $i(x) = (x,e,x)$.
We can define the product and inverse as follows: $(x,g,y)(y,h,z) = (x,gh,z)$ and
$(x,g,y)^{-1} = (y,g^{-1},x)$.
It is not hard to check that $G$ is a groupoid.
Groupoids that defined by group bundles in this way are said to be of \emph{standard type}.
The interesting thing is that every groupoid is actually of standard type.

\medskip

\begin{prop}
\label{gpd=gp-bdl}
There is a natural one to one correspondence between groupoids with unit space $X$ and
group bundles over the equivalence classes of equivalence relations on $X$.
\end{prop}
\noindent
\begin{prf}
We have already seen in the above how we can construct a groupoid with unit space $X$
from a group bundle over the equivalence classes of an equivalence relation on $X$ and this gives a
correspondence $\cal C$.
We give the inverse construction in the following.
Let $(H, X, \tar, \s, i)$ be a groupoid.
For any $x,y\in X$, we let
$$\ \!\!_xH_y\ :=\ \tar^{-1}(x)\cap \s^{-1}(y).$$
It is clear that $\ \!\!_xH_x$ is a group.
Moreover, if we take any $g\in \ \!\!_xH_y$, the map $h\mapsto g^{-1}hg$ is an isomorphism from
$\ \!\!_xH_x$ to $\ \!\!_yH_y$.
Furthermore, we have $\ \!\!_xH_x g = \ \!\!_xH_y = g \ \!\!_yH_y$.
We define an equivalence relation $S$ on $X$ by $(x, y)\in S$ if $\ \!\!_xH_y\neq \emptyset$.
For any $\xi\in X/S$, we fix a representative $x_\xi\in \xi$.
Then, we obtain an assignment of groups $\xi\mapsto \ \!\!_{x_\xi}H_{x_\xi}$.
It is clear that if the groupoid $H$ is defined by a group bundle
$\{G_\xi\}_{\xi\in X/R}$, then $S = R$
and $\ \!\!_{x_\xi}H_{x_\xi}\cong G_\xi$ for any $\xi \in X/R$.
Therefore, the correspondence $\cal C$ is injective.
On the other hand, we need to show that $\cal C$ is surjective.
It is required to show that starting from a groupoid $H$ with unit space $X$,
the canonical groupoid $G$ associate with the group bundle
$\{\ \!\!_{x_\xi}H_{x_\xi}\}_{\xi\in X/S}$ over the equivalence classes of the equivalence relation $S$
is isomorphic to $H$.
In fact, for any $\xi\in X/R$ and any $x\in \xi$, we fix an element
$l_x\in \ \!\!_{x_\xi}H_x$.
For any $g\in H$, let $x= \tar(g)$ and  $y = \s(g)$.
Since $(x, y)\in S$, they belong to an equivalence class $\xi$.
Now, we define
$$\Phi(g)\ :=\ (x, l_x \cdot g \cdot l_y^{-1}, y) \in G.$$
It is clear that $\Phi$ is a groupoid homomorphism.
Moreover, $\Phi$ is injective because $l_x \cdot g \cdot l_y^{-1} = l_x \cdot g' \cdot l_y^{-1}$
will imply that $g=g'$ and $\Phi$ is surjective since $\Phi(l_x^{-1} \cdot h \cdot l_y) =
(x,h,y)$ (for any $(x,h,y)\in G$). 
\end{prf}

\medskip

In the case when $H$ is only a small category, one can also define the equivalence relation $R$ as well as $_xH_y$ as in the above. 

\medskip

\begin{cor}
Let $\alpha$ and $\beta$ be actions of groups $G$ and $H$ on sets $X$ and $Y$ respectively.
Then $X\times_\alpha G \cong Y\times_\beta H$ as groupoids if and only if
there exists a bijection $\psi: X\rightarrow Y$ such that the restriction of $\psi$ on
every orbit $O$ of $\alpha$ is a bijection onto an orbit of $\beta$
whose stabilization group (a subgroup of $H$) is isomorphic to that of $O$ (a subgroup of $G$).
\end{cor}

\medskip

\begin{defn}
\label{def-act}
Let $G$ and $H$ be small categories. 
Suppose that $\varphi: H^{(0)} \rightarrow G^{(0)}$. 
We let 
$$G\times^\varphi H\ :=\ \{ (g,h)\in G\times H: \s(g) = \varphi(\tar(h)) = \varphi (\s(h)) \}.$$
A \emph{left action of $G$ on $H$ with respect to $\varphi$} is a map $(g,h)\mapsto \alpha_g(h)$ from $G\times^\varphi H$ to $H$ such that for any $(g',g)\in G^{(2)}$, $(h',h)\in H^{(2)}$ and $u\in H^{(0)}$ with $(g,h),(g,u), (g,h')\in G\times^\varphi H$, we have: 
\begin{enumerate}
\item[\rm (0).] $\s(\alpha_g(u)) = \tar(\alpha_g(u))$;
\item[\rm (I).] $\s(\alpha_g(\s(h))) = \s(\alpha_g(h))$ ;
\item[\rm (II).] $\tar(\alpha_g(\tar(h))) = \tar(\alpha_g(h))$; 
\item[\rm (III).] $\varphi(\s(\alpha_g(u))) = \tar(g)$;
\item[\rm (IV).] $\alpha_{\varphi(\tar(h))} (h) = h$;
\item[\rm (V).] $\alpha_{g'}(\alpha_{g}(h)) = \alpha_{g'g}(h)$;
\item[\rm (VI).] $\alpha_g(h' h) = \alpha_g(h')\alpha_g(h)$.
\end{enumerate}
For simplicity, we say that $(\varphi, \alpha)$ (or simply $\alpha$ if $\varphi$ is understood) is \emph{a left action of $G$ on $H$}. 
\end{defn}

\medskip

\begin{rem}
\label{rem-act}
{\rm 
(a) From {\rm (I)}, {\rm (II)} and {\rm (III)}, we know that 
\begin{equation}
\label{s-alp-h}
\varphi(\s(\alpha_g(h)))\ =\ \tar(g)\ =\ \varphi(\tar(\alpha_g(h))).
\end{equation} 
Therefore, 
$\varphi(\tar(\alpha_{g}(h))) = \tar (g) = \s(g')$ and $\varphi(\s(\alpha_{g}(h))) = \s(g')$
which show that the left hand side of {\rm (V)} makes sense. 
The right hand side of {\rm (V)} clearly makes sense because $\s(g'g) = \s(g)$. 

\smnoind
(b) It is clear that the left hand side of {\rm (VI)} makes sense as $\tar(h'h) = \tar(h')$ and $\s(h'h) = \s(h)$. 
On the other hand, since $\tar(\alpha_g(h)) = \tar(\alpha_g(\tar(h))) = \tar(\alpha_g(\s(h'))) = \s(\alpha_g(\s(h'))) = \s(\alpha_g(h'))$ (because of {\rm (0)}, {\rm (I)} and {\rm (II)}), the right hand side of {\rm (VI)} makes sense. 

\smnoind
(c) For any $y\in G^{(0)}$, we denote $X_y = \varphi^{-1}(y)$. 
For any $g\in G$, the left action $\alpha$ defines a natural transform
$$\alpha_g: H\!\mid\!_{X_{\s(g)}} \rightarrow H\!\mid\!_{X_{\tar(g)}}$$ 
(where $H\!\mid\!_Z$ is the full subgroupoid of $H$ with $H\!\mid\!_Z^{(0)} = Z$).  
If $y\in G^{(0)}$, then $\alpha_y$ is the identity map from $H\!\mid\!_{X_y}$ to $H\!\mid\!_{X_y}$. 
Moreover, if $G$ is a groupoid, then $\alpha_g$ is an ismorphism. 
}
\end{rem}

\medskip

\begin{lem}
\label{lem-act}
Suppose that $G$ is a category acting on another category $H$ through a left action $(\varphi, \alpha)$. 

\smnoind
(a) If $G_\varphi$ is the full subcategory of $G$ with $G_\varphi^{(0)} = \varphi(H^{(0)})$, then $$G\times^\varphi H\ =\ G_\varphi\times^\varphi H.$$
Consequently, we can always assume that $\varphi$ is surjective.

\smnoind
(b) If $G$ is a groupoid and $(g,u) \in G\times^\varphi H$ with $u\in H^{(0)}$, then $\alpha_g(u)\in H^{(0)}$. 
In this case, $\alpha$ induces an action $\alpha^{(0)}$ of $G$ on $H^{(0)}$. 

\smnoind
(c) Suppose that both $G$ and $H$ are groupoids and $(g,h)\in G\times^\varphi H$. 
Then 
\begin{equation}
\label{alp-inv}
\alpha_g(h)^{-1}\ =\ \alpha_g(h^{-1}). 
\end{equation}
\end{lem}
\begin{prf}
(a) For any $(g,h)\in G\times^\varphi H$, we have $\s(g) = \varphi (\s(h))\in G_\varphi^{(0)}$ and 
$\tar(g) = \varphi (\tar(\alpha_g(h)))\in G_\varphi^{(0)}$ (by Remark \ref{rem-act}(a)) and the equality follows.  

\smnoind
(b) If $v = \s(\alpha_g(u)) = \tar (\alpha_g(u))$, then 
$$\alpha_{g^{-1}}(\alpha_g(u))\ =\ \alpha_{g^{-1}}(\alpha_g(u)v)\ =\ \alpha_{\s(g)}(u)\alpha_{g^{-1}}(v)\ =\ u\alpha_{g^{-1}}(v)\ =\  \alpha_{g^{-1}}(v)$$ 
and so $\alpha_g(u) =v$. 

\smnoind
(c) It is clear that $(g,h^{-1})\in G\times^\varphi H$. 
We have by (II) and (IV), 
$$\alpha_g(h)\alpha_g(h^{-1})\ =\ \alpha_g(\tar(h))\ =\ \tar(\alpha_g(h))$$ 
and similarly, $\alpha_g(h^{-1})\alpha_g(h) = \s(\alpha_g(h))$. 
\end{prf}

\medskip

\begin{rem}
\label{gpd-act}
{\rm 
Suppose that $G$ is a groupoid. 
By Lemma \ref{lem-act}(b), one can replace Conditions {\rm (0) - (III)} with the following three conditions: 
\begin{enumerate}
\item[\rm (I').] $\alpha_g(\s(h)) = \s(\alpha_g(h))$;
\item[\rm (II').] $\alpha_g(\tar(h)) = \tar(\alpha_g(h))$;
\item[\rm (III').] $\varphi(\alpha_g(u)) = \tar(g)$.
\end{enumerate}
Note that condition {\rm (I')} implies that $\alpha_g(H^{(0)}) \subseteq H^{(0)}$ for every $g\in G$.
}
\end{rem}

\medskip

\begin{prop}
\label{cr-pd}
Suppose that $G$ is a groupoid acting on a small category $H$ by a left action $(\varphi, \alpha)$ and 
$$H\times_\alpha G\ :=\ \{(h,g)\in H\times G: \tar(g) = \varphi(\s(h)) = \varphi (\tar(h))\}.$$ 
For any $(h,g)\in H\times_\alpha G$, we set 
$$\s(h,g)\ :=\ \alpha_{g^{-1}}(\s(h)) \qquad {\rm and} \qquad \tar(h,g)\ :=\ \tar(h)$$
(here, we identify $u\in H^{(0)}$ with its canonical image $(u, \varphi(u)) \in H\times_\alpha G$). 
Moreover, if $(h,g),(h',g')\in H\times_\alpha G$ satisfying $\s(\alpha_{g^{-1}}(h)) = \tar(h')$, we define
$$(h,g)(h',g')\ :=\ (h\alpha_g(h'), gg').$$
This turns $H\times_\alpha G$ into a small category. 
If, in addition, $H$ is a groupoid, then $H\times_\alpha G$ is also a groupoid with 
$$(h,g)^{-1} = (\alpha_{g^{-1}}(h^{-1}), g^{-1}).$$
\end{prop}
\begin{prf}
Since 
\begin{equation}
\label{gg'}
\tar(g')\ =\ \varphi(\tar(h'))\ =\ \varphi(\s(\alpha_{g^{-1}}(h)))\ =\ \tar (g^{-1})\ =\ \s(g)
\end{equation}
(by Equality (\ref{s-alp-h})), the product $gg'$ is valid. 
Secondly, as $(h',g')\in H\times_\alpha G$, Equation (\ref{gg'}) shows that $(g,h')\in G\times^\varphi H$. 
Furthermore, 
\begin{equation}
\label{algh'}
\tar(\alpha_g(h'))\ =\ \alpha_g(\tar(h'))\ =\ \alpha_g(\s(\alpha_{g^{-1}}(h)))\ =\ \s(h)
\end{equation}
by the hypothesis as well as (I'), (II') \& (IV) and so $(h,\alpha_{g}(h'))\in H^{(2)}$. 
Therefore, the product is well defined. 
It is not hard to see that the product is associative. 
Now, suppose that $H$ is a groupoid. 
Then clearly $(g^{-1},h^{-1})\in G\times^\varphi H$, and by Equation (\ref{s-alp-h}), 
$$\varphi(\s(\alpha_{g^{-1}}(h^{-1})))\ =\ \tar(g^{-1})\ =\ \varphi(\tar(\alpha_{g^{-1}}(h^{-1}))).$$
Thus, $(\alpha_{g^{-1}}(h^{-1}), g^{-1})\in H\times_\alpha G$. 
Moreover, by (I') and (V), 
$$\s(\alpha_{g^{-1}}(h^{-1}), g^{-1})\ =\ \alpha_g(\s(\alpha_{g^{-1}}(h^{-1})))\ =\ \tar(h)\ =\ \tar(h,g)$$ 
and 
$$\tar(\alpha_{g^{-1}}(h^{-1}), g^{-1})\ =\ \tar(\alpha_{g^{-1}}(h^{-1}))\ =\ \alpha_{g^{-1}}(\tar(h^{-1}))\ =\ \s(h,g).$$ 
Now, it is easy to check that $(\alpha_{g^{-1}}(h^{-1}), g^{-1})$ is the inverse of $(h,g)$. 
\end{prf}

\medskip

\begin{rem}
{\rm
(a) Suppose that $\bar H = \{\bar h: h\in H\}$ is the opposite category of $H$ (i.e. $\s(\bar h) = \tar (h)$, $\tar (\bar h) = \s(h)$, $\bar H^{(2)} = \{(\bar h,\bar h'): (h',h)\in H^{(2)}\}$ and $\bar h \bar h' = \overline{h'h}$). 
Since $\bar H^{(0)} = H^{(0)}$, it is not hard to check that a left action $(\varphi, \alpha)$ of $G$ on $H$ induces a left action $(\bar \alpha, \varphi)$ on $\bar H$ given by $\bar \alpha_g(\bar h) = \overline{\alpha_g(h)}$. 

\smnoind
(b) One can also turn $G\times^\varphi H$ into a small category with the following source map, target map and product:
$$\s(g,h)\ :=\ \s(h), \quad \tar(g,h)\ :=\ \tar(\alpha_g(h)) \quad {\rm and} \quad 
(g',h')(g,h)\ :=\ (g'g, \alpha_{g^{-1}}(h')h)$$
(if $\s(h') = \tar(\alpha_g(h))$). 
If $\bar H$ and $\bar \alpha$ is as in part (a), then it is not hard to check that the map $\Psi: \bar H \times_{\bar \alpha} G \rightarrow G\times^\varphi H$ given by $\Psi(\bar h, g) = (g^{-1}, h)$ is an isomorphism of the two categories. 
Consequently, if $H$ is a groupoid, the map $\Psi:H\times_\alpha G \rightarrow G\times^\varphi H$ defined by $\Psi(h,g) = (g^{-1},h^{-1})$ is a groupoid isomorphism. 
}
\end{rem}

\medskip

We require that $G$ is a groupoid in Proposition \ref{cr-pd} because we need to define $\s(h,g)$. 
There is a situation when $\s(h,g)$ can be determined without the existence of $g^{-1}$. 
Note that since $G$ is groupoid, $\varphi(\alpha_{g^{-1}}(\s(h))) = \tar(g^{-1}) = \s(g)$ (by (III')) and one has $\alpha_{g^{-1}}(\s(h)) = \varphi^{-1}(\s(g))$ if $\varphi$ is injective. 
However, in this case, one can even assume that $\varphi$ is actually bijective
(because of Lemma \ref{lem-act}(a)), and one can identify $H^{(0)}$ with $G^{(0)}$. 

\medskip

The following proposition follows from nearly the same argument as that of Proposition \ref{cr-pd} except that we need to replace equalities (\ref{algh'}) by the following: 
$$\tar(\alpha_g(h'))\ =\ \tar(\alpha_g(\tar(h')))\ =\ \s(\alpha_g(\tar(h')))\ =\ \tar(g)\ =\ \s(h).$$

\medskip

\begin{prop}
\label{cr-pd-id}
Let $G$ and $H$ be small categories such that $H^{(0)} = G^{(0)}$. 
Suppose that $(\id, \alpha)$ a left action of $G$ on $H$. 
For any $(h,g)\in H\times_\alpha G := \{(h,g)\in H\times G: \tar(g) = \s(h) = \tar(h)\}$, we set 
$$\s(h,g)\ :=\ (\s(g),\s(g)) \qquad {\rm and} \qquad \tar(h,g)\ :=\ (\tar(h),\tar(h)).$$
Moreover, if $(h,g),(h',g')\in H\times_\alpha G$ satisfying $\s(g) = \tar(h')$ we define
$$(h,g)(h',g')\ :=\ (h\alpha_g(h'), gg').$$
This turns $H\times_\alpha G$ into a small category. 
\end{prop}

\medskip

The category $H \times_\alpha G$ is called the \emph{semi-direct product of $H$ and $G$ under the left action $(\varphi, \alpha)$} (when $G$ is a groupoid or when $H^{(0)} = G^{(0)}$). 

\medskip

Now, let $G$ be a category acting on a category $H$ through a left action $(\varphi, \alpha)$ and $R$ be the canonical equivalence relation on $G^{(0)}$ defined by $G$. 
If $G^{(0)} = H^{(0)}$ and $\varphi = \id$, then $_v(H\times_\alpha G)_u = \ \!\!_vH_v\times \ \!\!_vG_u$ and the composition is given by the map $_vG_u\times \ \!\!_uH_u \rightarrow \ \!\!_vH_v$ ($u,v\in H^{(0)}$).  
Thus, in order to construct $H\times_\alpha G$, one needs only to know the semi-group bundles $_uH_u$ ($u\in H^{(0)}$). 
One can also use the idea from this decomposition to construct a ``restricted semi-direct product''. 

\medskip

For any $u,v\in H^{(0)}$ with $\varphi(u) R \varphi(v)$, we define 
$$_u\ti G_v\ :=\ \{ g\in \ \!\!_{\varphi(u)} G _{\varphi(v)}: u =\alpha_g(v) \}$$
(which could be empty). 
For any $g\in \ \!\!_u\ti G_v$ and $h\in \ \!\!_v\ti G_w$, we have $gh \in \ \!\!_{\varphi(u)} G _{\varphi(w)}$ and $\alpha_{gh}(w) = \alpha_g(\alpha_h(w)) = u$ and thus, $gh\in \ \!\!_u\ti G_w$.  

\medskip

Let $\ti G$ be category with object $H^{(0)}$ and with morphisms from $u\in H^{(0)}$ to $v\in H^{(0)}$ being $_u\ti G_v$ (note that the identity morphism in $_u\ti G_u$ is $\varphi(u)\in \ \!\!_u\ti G_u$). 

\medskip

Consider $j$ be the canonical natural transform from $\ti G$ to $G$ that send an object $u$ to $\varphi (u)$ and send a morphism $g\in \ \!\!_u\ti G_v$ to $g\in \ \!\!_{\varphi(u)} G _{\varphi(v)}$. 
Therefore, for any $g\in \ \!\!_u\ti G_v$, we have $\s(j(g)) = \varphi(\ti \s(g))$ and $\tar(j(g)) = \varphi(\ti \tar(g))$ (where $\ti \s$ and $\ti \tar$ are the source and the target maps in $\ti G$). 
It is clear that for any $(g,h)\in \ti G\times^\id H$, we have $(j(g), h) \in G\times^\varphi H$ and we can define $\ti \alpha_g(h) = \alpha_{j(g)}(h)$. 
It is not hard to see that $(\id, \ti \alpha)$ is an action of $\ti G$ on $H$. 
By Proposition \ref{cr-pd-id}, one can define $H \times_{\ti \alpha} \ti G$ which is called the \emph{restricted semi-direct product of $H$ and $G$ under the action $\alpha$}.

\medskip

\begin{rem}
\label{re-cr-pd}
{\rm 
(a) If $G$ is a groupoid, then for any $g\in \ \!\!_u\ti G_v$, we also have $g^{-1}\in \ \!\!_{\varphi(v)} G _{\varphi(u)}$ and $v = \alpha_{g^{-1}} (u)$ which shows that $g^{-1} \in \ \!\!_v\ti G_u$ and so, $\ti G$ is also a groupoid. 
Note that $(h, j(g)) \in H\times_\alpha G$ if $(h, g) \in H\times_{\ti \alpha} \ti G$ and the map $\Psi$ that send $(h,g)$ to $(h,j(g))$ is an injective natural transform that is an identity on the space of object $H^{(0)}$. 
In general, $\Psi$ is not surjective because for any $(h,g)\in H \times_{\ti \alpha} \ti G$, one has $\s(h) = \tar(h)$. 
That is why $H \times_{\ti \alpha} \ti G$ is called the restricted semi-direct product. 
Note that if $G$ and $H$ are groups, then $H \times_{\ti \alpha} \ti G = H \times_\alpha G$. 

\smnoind
(b) Suppose that $G$ is a groupoid and $H'$ is the semi-group bundle of $H$ (i.e. $_uH'_u = \ \!\!_uH_u$ and $_uH'_v = \emptyset$ if $u\neq v$). 
Then $(\varphi, \alpha)$ induces an action $(\varphi, \alpha')$ of $G$ on $H'$. 
Obviously, $H'\times_{\alpha'} G$ is a subcategory of $H \times_\alpha G$. 
Moreover, it is clear that $H'\times_{\alpha'} G$ contains the image of the natural transform $\Psi$ in part (a). 
Conversely, for any $(h,g)\in H'\times_{\alpha'} G$, we see that $\tar (g) = \varphi(\s(h))$ and $\s(h) = \tar(h)$. 
We let $u =\s(h)$ and $v = \alpha_{g^{-1}}(u)$. 
Then $\varphi(v) = \s(g)$ and so $g\in \ \!\!_u\ti G_v$. 
Thus, $\Psi(H\times_{\ti \alpha} \ti G) = H'\times_{\alpha'} G$. 
Consequently, if $H$ is a semi-group bundle, then $H \times_{\alpha} G \cong H\times_{\ti \alpha} \ti G$ and this applies, in particular, to the case when $H$ is a set. 
Now, if $G$ is only a category and $H$ is a semi-group bundle, then one can also define $H \times_{\alpha} G$ as $H\times_{\ti \alpha} \ti G$. 
}
\end{rem}

\medskip

\begin{eg}
\label{eg-act}
(a) Suppose that $G$ is a group and $H$ is the groupoid given by an equivalence relation $\sim$ on $X = H^{(0)}$. 
Then a left action of $G$ on $H$ is the same as an action $\beta$ of $G$ on $X$ that respects $\sim$ in the following sense: for any $u,v\in X$ and $g\in G$, 
$$u\sim v \qquad {\rm implies} \qquad \beta_g(u)\sim \beta_g(v).$$

\smnoind
(b) If both $G$ and $H$ are group (i.e. both $G^{(0)}$ and $H^{(0)}$ are singletons), then action of $G$ on $H$ in the above coincides with the usual definition of action of $G$ on $H$ by automorphisms and $H\times_\alpha G$ is the usual semi-direct product. 

\smnoind
(c) If $G$ is a groupoid and $H$ is a set $X$ together with the trivial groupoid structure (i.e. $H^{(2)} =  \{(x,x): x\in X\}$), then action of $G$ on $H$ in the above coincides with the usual definition of action of $G$ on the set $X$ and $X\times_\alpha G$ is the usual semi-direct product of $G$ and $X$. 

\smnoind
(d) If $G$ is a groupoid and $G = G^{(0)}$, then $H\times_\alpha G \cong \{ h\in H: \varphi(\s(h)) = \varphi(\tar(h))\}$ as categories. 
In particular, if $\varphi(H^{(0)})$ is a singleton, then $H\times_\alpha G \cong H$. 
On the other hand, if $\varphi$ is injective, then $H\times_\alpha G$ coincides with the semi-group bundle $\{ h\in H: \s(h) = \tar(h)\}$. 

\smnoind
(e) Let $G$ be a groupoid, $H = G$ and $\varphi$ be the identity map. 
We define an action of $G$ on itself by inner automorphism, i.e. $\alpha_g(h) = ghg^{-1}$ (for any $(g,h)\in G\times^\id G$). 
It is not hard to check that the conditions in Definition \ref{def-act} are satisfied. 
Moreover, as in the case of group, there is a canonical groupoid homomorphism $\Psi$ from $G\times_\alpha G$ to $G$ that sends $(h,g)$ to $hg^{-1}$. 
However, unlike the group case, the ``kernel of $\Psi$'' equals 
$$\{ (h,g)\in G\times_\alpha G: \Psi(h,g) \in G^{(0)} \}\ =\ \{ (g,g^{-1})\in G\times G: \s(g) = \tar(g) \} \ \cong\ \{g\in G: \s(g) = \tar(g) \}$$ and is not isomorphic to the whole of $G$ (unless $G^{(0)}$ is a singleton). 

\smnoind
(f) Let $X$ be a set of Hilbert spaces, $H$ is the category of bounded linear maps between elements in $X$ and $G$ is the category of isometries from one element of $X$ to another one. 
We define a left action $\alpha$ of $G$ on $H$ as follows. 
If $\Psi: K_1 \rightarrow K_2$ is an isometry and $T: K_1\rightarrow K_1$ is a bounded linear map, then $\alpha_\Psi(T) = \Psi\circ T \circ \Psi^*$. 
It is not hard to see that $(\id, \alpha)$ satisfy the conditions in Definition \ref{def-act}. 
For any $(T,\Psi) \in H\times_\alpha G$, we let $\Phi(T,\Psi) = T\circ \Psi^*$. 
Then it is not hard to see that $\Phi$ is a natural transform from $H\times_\alpha G$ to $H$. 
In general, $\Phi$ is surjective but not injective. 
\end{eg}

\medskip

\noindent
Chern Institute of Mathematics and LPMC,
Nankai University,
Tianjin 300071,
China.

\smallskip\noindent
E-mail address: ckng@nankai.edu.cn
\end{document}